\documentclass[a4paper, 11pt, p1, authoryear]{article}

\usepackage{tabularx}
\usepackage[centerlast,footnotesize,width=.975\linewidth]{caption}
\usepackage[labelformat=simple]{subcaption}

\usepackage{hyperref}
\usepackage{float}
\usepackage{booktabs}
\usepackage{amsmath, amsfonts, amssymb}

\usepackage[margin=2cm, nofoot]{geometry}
\usepackage{fancyhdr}
\usepackage{color,soul}
\usepackage[toc,page]{appendix}
\usepackage[export]{adjustbox}
\usepackage{listings}
\usepackage[boxed,linesnumbered]{algorithm2e}
\SetKwRepeat{Do}{do}{while}
\definecolor{mygreen}{rgb}{0,0.6,0}
\definecolor{mymauve}{rgb}{0.58,0,0.82}

\usepackage{amsthm}
\usepackage[nameinlink]{cleveref}

\theoremstyle{definition}
\newtheorem{definition}{Definition}[section]
\theoremstyle{remark}
\newtheorem{remark}{Remark}

\graphicspath{ {plots/}}

\usepackage[utf8]{inputenc}  
\DeclareMathOperator{\mydiv}{div}

\usepackage{authblk}
\usepackage{natbib}

\begin{document}
\title{Layer-adapted meshes for singularly perturbed problems via mesh
  partial differential   equations   and \emph{a posteriori} information}
\author[1]{R\'ois\'in Hill}%
\author[2]{Niall Madden}
\affil[1]{Department of Mathematics \& Statistics, University of
  Limerick, Ireland. Corresponding author. Email:~\texttt{roisin.hill@ul.ie}. }
\affil[2]{School of Mathematical and Statistical Sciences,
  University of Galway, Ireland. Email:~\texttt{Niall.Madden@UniversityofGalway.ie}}
\date{\today}
\maketitle
\thispagestyle{empty}  

\begin{abstract}
We propose a new method for the construction of layer-adapted
meshes for singularly perturbed differential equations (SPDEs), based on
mesh partial differential equations (MPDEs) that incorporate  \emph{a
  posteriori} solution information.
There are numerous studies on the development of parameter robust
numerical methods for SPDEs that depend on the layer-adapted mesh of
Bakhvalov. In~\citep{HiMa2021}, a novel MPDE-based approach for
constructing a generalisation of these meshes was proposed.
Like with most layer-adapted mesh methods, the algorithms in that
article depended on
detailed derivations of \emph{a priori} bounds on the SPDE's solution
and its derivatives.
In this work we extend that approach  so that it
instead uses  \emph{a  posteriori} computed estimates of the solution.
We present detailed algorithms for the efficient implementation of the
method, and numerical results for the robust solution of two-parameter
reaction-convection-diffusion problems, in one and two dimensions.
We also provide full FEniCS code for a one-dimensional example.
\end{abstract}

\section{Introduction}
\label{sec: intro MPDE dx}
This article is concerned with a new approach to generating
layer-adapted meshes for singularly perturbed differential equations (SPDEs).
The core ideas is to use a new formulation for the classic fitted
meshes of Bakhvalov~\citep{Bakh1969} proposed in~\citep{HiMa2021},
but extended to use \emph{a  posteriori} computed quantities, rather
than the usual  \emph{a priori} information usually used to construct
these meshes.

For exposition, we will focus on the numerical solution of
two-parameter reaction-convection-diffusion equations of the form
\begin{equation}%
  \label{eqn: model RCD}
  -\varepsilon \Delta u(x) + \mu\mathbf{b}\cdot\nabla u(x) + r(x) u(x) =  f(x)\quad \text{for }x\in 
  \Omega^{d},\quad\text{ with } u\rvert_{\partial \Omega} = 0,
\end{equation}
with $d=1,2$. We make rather standard assumptions on the other problem
data; specifically that $\mathbf{b}$, $r$, and $f$ are given smooth
functions, and that $2r >\mu\mydiv\mathbf{b}$. 
Equation~\eqref{eqn: model RCD} features a pair of positive
parameters, $\varepsilon$ and $\mu$, which may be arbitrarily small, making this a
\emph{singularly perturbed} problem.
Typically, solutions to \ref{eqn: model RCD} exhibit layers,
the location and width of which can be challenging to determine
\emph{a~priori} (especially when $d=2$), making it particularly interesting for exhibiting the features of our proposed
method.  The
proposed method automatically determines these quantities, and
constructs an appropriate mesh by solving a suitable  mesh partial differential
equation (MPDE). This is done in practice by alternating between
solving the SPDE and the MPDE on a nested sequence of grids using
standard Galerkin  finite element methods.

SPDEs, such as \eqref{eqn: model RCD},   are of interest to mathematical
modellers, since they can be applied to describe a wide range of
physical phenomena.  Their numerical solution is of
significant interest in numerical analysis, and great efforts have
been made to devise (and analyse) methods which can solve such
problems accurately, and resolve any layers present. A detailed
overview of the field (as it was in 2008) in given in~\citep{RoSt2008}; see also
\citep{Roos2022} for a more recent view of advances and challenges.

One of the challenges in the numerical solution of singularly
perturbed problems is the development of methods for which a
meaningful error bound can be established that is independent of the
perturbation parameter(s), and which ensure any layers present are
resolved; the monograph of \citep{MiOR2012} provide detailed motivation
for this and presents methods that enjoy these properties,
for a wide class of problems.  These methods are mainly based on the
famous piecewise uniform Shishkin mesh~\citep{ShSh2009}. We refer to
\citep{Lins2010} for a more general treatment, which includes analyses
for other meshes, including the graded Bakhvalov
mesh~\citep{Bakh1969}.

The meshes mentioned above are constructed based on \emph{a priori}
information on the solution and its derivatives. The approach that we
present is closer in philosophy to \emph{a posteriori} adaptive
algorithms, of which there are many in the literature; notable
examples include the now-classic work of~\cite{BeMa2000}
and~\cite{KoSt2001}.
We also mention \cite{SiSh2013,ShMo2022},
which are closer in style to this article, since they use moving mesh
methods, as well as the work on reaction-convection-diffusion problems
of~\cite{WuZh13}.
However,
the method that we propose is distinguished in that our goal is to
automatically reconstruct the mesh density function of a Bakhvalov
mesh, rather than by adapting the mesh directly.

We emphasise that our goal is not to present an algorithm for
\eqref{eqn: model RCD}, \emph{per se}, but to use it as a test case
for testing our approach. Furthermore, since the details of the
construction of the usual Bakhvalov mesh may differ
substantially from when $d=1$ to $d=2$, we defer a detailed
description (and review of the literature) to \S\ref{sec: 1D SPDE} and
\S\ref{sec: 2D SPDE}, respectively.

The rest of this article is organised as follows. In the next section,
we summarise some notation used throughout.
We then turn our attention to one-dimensional  versions of \eqref{eqn:
  model RCD} in \S\ref{sec: 1D problems}. The background and some key
references are discussed in \S\ref{sec: 1D SPDE}
In \S\ref{sec: MPDE
   dx} we present the MPDE that we use to generate the meshes on which
 to solve the SPDE.
 We present the algorithm
 to generate these meshes and detail how our method is implemented in
 \S\ref{sec: Algorithm for MMPDE dx}. In \S\ref{sec:1D numerics} we present
 the results of numerical experiments which 
 verify the accuracy and efficiency of the method.
 
 In \S\ref{sec: 2D problem} we extend the approach to two-dimensional
  reaction-convection-diffusion problems, where, in addition to the
  relationship between the two parameters, the direction of the flow
  influences the nature and location of layers.
  We present the MPDE formulations in \S\ref{sec: 2D MPDE dx}, and
  the implementation in \S\ref{sec: 2D algorithm dx V2}. Again,
  validating numerical results are presented in \S\ref{sec: 2D
    numerical results}.
We use a standard Galerkin $\mathcal{P}_1$ finite element method
to solve both the MPDEs and SPDEs. All results presented have been computed using
FEniCS~\citep{LoMa2012}, with a full working example presented in
\autoref{app: code 1D}. We present 
pseudo-code the for two-dimensional problem in \autoref{app: 2D algorithm}.

\subsection*{Notation and definitions}
We use $\Omega$ to denote an open, bounded subset of $\mathbb{R}^{d}$,
$\overline{\Omega}$ denotes its closure, and $\partial \Omega$ its
boundary.  Usually, it is the domain on which a (physical) SPDE
is posed. We denote the computational domain (on which the MPDE is posed) as $\Omega^{[c]}$. 

Where necessary, we use a superscript to indicate a
differential operator applied on the computational domain. In
particular, we define
\begin{equation*}
  \nabla^{[c]} = \left(\frac{\partial}{\partial \xi_1}, \frac{\partial}{\partial \xi_2}\right),\quad\text{for } (\xi_1, \xi_2)\in \Omega^{[c]}.
\end{equation*}

We use $\omega_h$ to denote a mesh that discretizes $\overline{\Omega}$, and $\omega_h^{[c]}$ as a mesh on $\Omega^{[c]}$. 
That is, it  denotes a partition, for the purposes of finite
element discretization,  into simplices of $\overline{\Omega}$, i.e.,
intervals in one dimension and triangles in two dimensions.
In addition, $\omega_h^{[c]}$ denotes a partition of a
computational domain, $\overline{\Omega}^{[c]}$.
A member of a sequence of partitions of $\overline{\Omega}$ or of
$\overline{\Omega}^{[c]}$ is denoted $\omega_h^{[i]}$ or
$\omega_h^{[c,i]}$.

\begin{definition}[The one-dimensional equidistribution principle]%
  \label{Def: Equidistribution principle}
  Let  $\rho$ := $\overline{\Omega}\rightarrow \mathbb{R}_{>0}$ be a
  strictly positive function known as the \emph{mesh density
    function}.
  We say that the mesh $\omega_h := \{a = x_0<x_1<\cdots<x_{N-1}<x_N =b\}$  \emph{equidistributes} $\rho$, if
  \begin{equation*}
    \int_{x_{i-1}}^{x_{i}} \rho(x)dx = \frac{1}{N}\int_{\overline{\Omega}} \rho(x)dx\quad \text{for }i = 1,\dots, N.
  \end{equation*}
\end{definition} 

\begin{definition}[Mesh generating function]%
  \label{Def: mesh generating func}
  A mesh generating function is a strictly monotonic bijective
  function $\varphi : \overline{\Omega}^{[c]}:= [0,1] \rightarrow
  \overline{\Omega}:=[a,b]$ that maps a uniform mesh with mesh points
  $\xi_i = i/N$, for $i = 0,1,\dots,N$, to a (possibly non-uniform)
  mesh with mesh points $ x_i = \varphi(i/N)$, for $i = 0,1,\dots,N$,
  with $\varphi(0) = a $ and $\varphi(1) = b$.
\end{definition}

\section{One-dimensional problems}%
\label{sec: 1D problems}
\subsection{A one-dimensional SPDE}\label{sec: 1D SPDE}
In this  section, we focus on the generation of meshes for solving the
one-dimensional reaction-convection-diffusion
problem
\begin{multline}%
  \label{eqn:1D RCD}
  -\varepsilon u''(x) + \mu b(x) u'(x) + r(x) u(x) =  f(x)\\ \text{for }
  x \in (0,1), \quad\text{ and } u(0)=u(1)=0.
\end{multline}
When $\varepsilon$ is small, and $\mu$ is
$\mathcal{O}(1)$ and positive, and assuming that $f$ does not vanish
at either boundary, a layer of width $\mathcal{O}(\varepsilon)$
will typically form on the right of the domain. If $\mu$ is negative, the layer
would be manifested near  the left boundary. However, if $\mu \ll 1$,
which makes this a so-called ``two-parameter'' problem,
then the situation is more complicated, and there may be layers at
both boundaries, whose widths depend on the relative magnitude of
$\varepsilon$ and $\mu$.

In spite of their apparent simplicity, one-dimensional linear problems
such as~\eqref{eqn:1D RCD} are widely studied (since, at least, the
work of~\citep{OMall1967}). In the numerical
analysis literature, progress was made in the early 2000s (see, e.g.,
\citep{RoUz2003,GrOR2006}). Research into these problems continues;
see, e.g.,  the analyses of a
discontinuous Galerkin method on \emph{a priori} layer-adapted
meshes~\citep{SiNa2020}, and an investigation of 
uniform convergence and supercloseness for the $\mathcal{P}_1$-FEM
solution on a graded meshes~\citep{ZhLv2022}.

Of particular interest to us  is the analysis  of a continuous Galerkin FEM applied  on a
Bakhvalov mesh~\citep{BrZa2016}. The level of detail in that paper
demonstrates the complexity in even constructing a suitable mesh for
this problem. The reason for this complexity is due to the interplay
between the values of $\varepsilon$ and $\mu$ which
determines the location and width of layers. This complicates the specification of \emph{a priori}
layer-adapted meshes, but not for the approach we propose, 
which automatically generates layer-adapted meshes for this
problem without  prior knowledge of the relationship between the
two-parameters.

We present the general approach (for any one-dimensional problem
featuring boundary layers) in \S\ref{sec: MPDE dx}, with more specific details of
the algorithm given in \S\ref{sec: Algorithm for MMPDE dx}.
In \S\ref{sec:1D numerics} we show how the approach
produces suitable meshes for a challenging two-parameter
problem with layers of different widths.

\subsection{A one-dimensional MPDE for Bakhvalov meshes}%
\label{sec: MPDE dx} 

In the literature, there are two distinct (but equivalent) approaches
to constructing meshes of Bakhvalov type. The original approach
depends on solving a certain non-linear problem (see, e.g.,
\citep[\S2.1.1]{Lins2010} and \citep[\S2.2]{HiMa2021}).
As noted in \citep[\S6.3]{Lins2010}, these meshes can also be constructed by
equidistributing certain monitor functions.
For a problem such as~\eqref{eqn:1D RCD}, where the solution may have a
layer near each of the boundaries, one equidistributes a mesh density function
of the form 
\begin{equation}\label{eq:rho 1DRD}
  \rho(x) = \max\bigg{\{}1,
  K \lambda_0\exp(-\lambda_0 x \sigma^{-1}) +
  K\lambda_1\exp(-\lambda_1 (1-x) \sigma^{-1})\bigg{\}},
\end{equation}
where $\lambda_0$ and $\lambda_1$ are related to the width of the
layers, $\sigma$ is determined by the formal order of the scheme, and
$0<K<1$ is chosen to control the proportion of the mesh points in the
layer regions.

It known that the mesh generating function for a mesh on
$\overline\Omega:=[a,b]$ that
equidistributes an arbitrary mesh generating function, $\rho$, 
can be expressed as the solution to the mesh PDE
\begin{equation}%
  \label{eqn: 1DMPDE}
  (\rho(x)x_{\xi}(\xi))_{\xi} = 0,\text{ for } \xi \in \Omega^{[c]}:= (0,1), \text{ with } x(0) = a \text{ and } x(1) = b;
\end{equation}
see \citep{HiMa2021} for details.
In the approach in that paper,~\ref{eqn: 1DMPDE} is solved
numerically with a $\mathcal{P}_1$-FEM, to obtain a mesh
generating function.
For the  classical construction of a Bakhvalov mesh,
$\rho$ depends on \emph{a priori} information concerning
the solution of the physical differential equation; in
\eqref{eq:rho 1DRD} it can be thought of as 
representing point-wise bounds for $|u'(x)|$.  

We now present an alternative approach, where $\rho=\rho(x,u_h)$ is
based on derivatives of numerical solutions 
to the SPDEs.
The resulting algorithm anticipates that there may be boundary
layers present, but does not require
\emph{a priori} knowledge concerning the boundaries at which layers may be present, or of the width of those layers.
The MPDE \eqref{eqn: 1DMPDE} is reformulated as
\begin{equation}%
\label{eqn: 1DMPDE dx}	
  (\rho(x, u_h)x_{\xi}(\xi))_{\xi} = 0,\text{ for } \xi \in \Omega^{[c]},
  \text{ with } x(0) = 0 \text{ and } x(1) = 1, 
\end{equation}
to emphasise that  now $\rho$ is dependent on the numerical solution
to the physical differential equation.
For~\eqref{eqn:1D RCD}, we take
\begin{subequations}
  \label{eqn: rho dx}
  \begin{equation}%
  \rho(x,u_h) = \max\bigg\{1,
  K\big(\upsilon_0\exp(-\upsilon_0 x\sigma^{-1}\big)\\
  +\upsilon_1\exp(-\upsilon_1(1-x)\sigma^{-1})\big)\bigg\},
\end{equation}
with
\begin{equation}\label{eq:upsilon_i}
  \upsilon_0(x, u_h) = \frac{\big|(u_h)^{}_{x^+}(x_0)\big|}
  {\max\big\{1, |f(x_0)|\big\}} \quad\text{ and }\quad
  \upsilon_1(x, u_h) = \frac{\big|(u_h)^{}_{x^-}(x_N)\big|}
  {\max\big\{1, |f(x_N)|\big\}},
\end{equation}
\end{subequations}
where  $f$ is the right-hand side of the SPDE.
It should be noted that this formulation is independent of the number
of boundary layers present in the solution, and their location(s) and
width(s). For example, if the solution possesses just one
layer near, which is near $x=1$,
then $(u_{h})_{x^-}(x_N)\gg 1$ leading to a graded mesh near that
boundary; elsewhere one will have $\rho(x,u_h)=1$, giving a uniform
coarse mesh over the rest of the domain.

\begin{remark}
  The term $\max\{1, |f(x_0)|\}$ in \eqref{eq:upsilon_i} handles the
  case where, for example, $f(0)=0$. Note that, in such an instance,
  the solution to \eqref{eqn:1D RCD} would not exhibit a strong
  boundary layer near $x=0$, in the sense that $u'(0)$ would be
  bounded independently of $\varepsilon$. Higher derivatives of $u$ at $x=0$
  would not be bounded, so some minor modifications of the approach
  may be need if using, for example, high-order finite elements.
\end{remark}

\subsection{Algorithm and implementation}%
\label{sec: Algorithm for MMPDE dx}

In \Cref{Alg: 1D MPDE dx} we outline a mechanism for solving the
nonlinear equation~\eqref{eqn: 1DMPDE dx} in an efficient manner.
At its core is a fixed-point iterative method, but it is accelerated
by using $h$-refinement~\citep{HiMa2021}. The inner iteration involves solving both
the SPDE and MPDE, since the latter requires  accurate estimates of the
derivatives of the solution to the SPDE.

The process begins by solving the SPDE on a uniform mesh with 16
elements (with fewer, a layer cannot be detected even for relatively
large values of $\varepsilon$). This is used to compute $\rho$ in \eqref{eqn:
  rho dx}, and thence a linearisation of~\eqref{eqn: 1DMPDE}. This
yields a new mesh on which the SPDE can be solved, and this process is
iterated until
\begin{equation*}%
  \label{eqn: change in dx}
  \Delta^{[k]} :=
  \Bigg{\lvert}\frac{(u_{h}^{[k]})_{x^+}(x_0)-(u_{h}^{[k-1]})_{x^+}(x_0)}
  {(u_{h}^{[k]})_{x^+}(x_0)}\Bigg{\rvert}
  +
  \Bigg{\lvert}\frac{(u_{h}^{[k]})_{x^-}(x_N)-(u_{h}^{[k-1]})_{x^-}(x_N)}
  {(u_{h}^{[k]})_{x^-}(x_N)}\Bigg{\rvert},
\end{equation*}
is less than some chosen tolerance. Then the mesh is uniformly refined,
and the iterative process repeated, until the computational mesh has
the desired number of intervals.
 
\begin{algorithm}[H]
	\SetKwInOut{Inputs}{Inputs}
	\caption{Generate a layer-adapted mesh using an MPDE, $h$-refinement and \emph{a posteriori} information.}
	\label{Alg: 1D MPDE dx}
	\Inputs{$N$, the number of intervals in the final mesh; $\rho$,
          a mesh density function; $\mathrm{TOL}$.}
	Set $\omega_h^{[c;0]} := \{\xi_0,\xi_1, \dots,\xi_{16}\}$ to be the uniform mesh on $\overline{\Omega}^{[c]}$ with $16$ intervals\;
	Set $x(\xi) = \xi$ for $\xi \in \overline{\Omega}^{[c]}$\;
	$r \leftarrow x$;  \quad	$k \leftarrow 0$\;			
	\For  { $i$ \emph{in} $0$\emph{:}$(\log_2(N){-}4)$}{
		\Do{ $\Delta^{[k]}  >  \mathrm{TOL}$}{
			Set $\omega_h^{[k]} := x(\omega_h^{[c;i]})$ to be the adapted mesh
			on $\overline{\Omega}$ at step $k$\;
			$u_h^{[k]} \leftarrow \mathcal{P}_1$-FEM solution to the SPDE on $\omega_h^{[k]}$\;
			Calculate $\Delta^{[k]}$\;
			Calculate $\upsilon_l(r, u_h^{[k]})$ for $l = 0,1$, and then $\rho(r, u_h^{[k]})$\;		
			set $x$ to be the $\mathcal{P}_1$-FEM
			solution, on $\omega_h^{[c;i]}$, to
			\begin{equation*}%
				\label{eqn: Art. MPDE2 dx}
				(\rho(r, u_h^{[k]}) x'(\xi))' = 0,\text{ for } \xi \in \Omega^{[c]} , x(0) = a \text{ and } x(1) = b;
			\end{equation*}
			$r \leftarrow  x$; \quad $k \leftarrow k{+}1$\;		
		}
		\If {$|\omega_h^{[c;i]}|<N{+}1$} {
			$\omega_h^{[c;i+1]}\leftarrow$ uniform $h$-refinement of
			$\omega_h^{[c;i]}$\;
			$r \leftarrow$ $r$ interpolated onto $\omega_h^{[c;i+1]}$\;		
	}}	
	Set $\omega_h := x(\omega_h^{[c;i]})$ to be the (final) adapted mesh
	on $\overline{\Omega}$.\;
\end{algorithm}

\subsection{Numerical results}
\label{sec:1D numerics}
\begin{figure}[H]
	\centering
	\begin{subfigure}{0.3\linewidth}
	\centering	
	\includegraphics[width=0.95\linewidth]{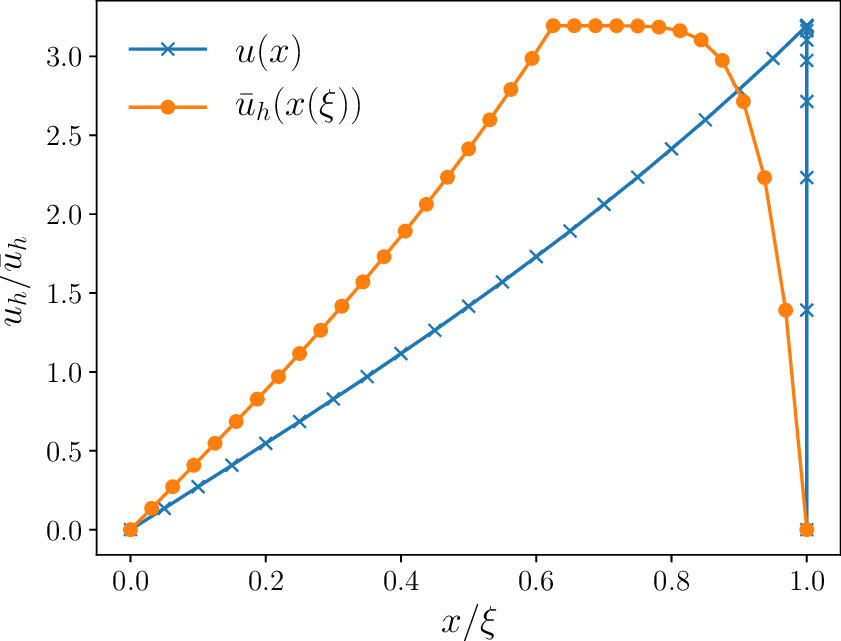}
	\caption{$\varepsilon=10^{-8}$ and $\mu = 1$}
	\label{Fig: 1DRCD mu1}
\end{subfigure}
\begin{subfigure}{0.3\linewidth}
	\centering
        \includegraphics[width=0.95\linewidth]{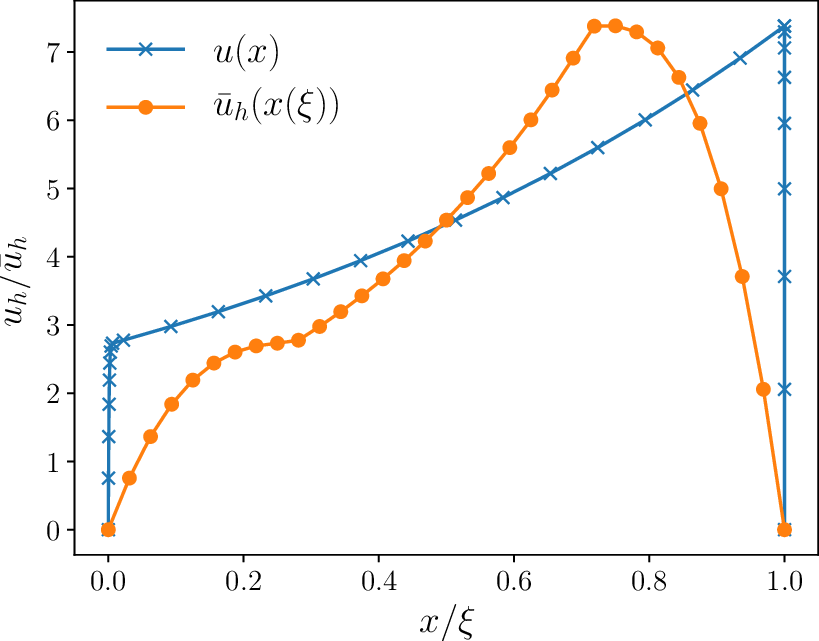}
	\caption{$\varepsilon=10^{-8}$ and $\mu = 10^{-3}$}
	\label{Fig: 1DRCD mu=0.0001}
\end{subfigure}
\begin{subfigure}{0.3\linewidth}
	\centering
	\includegraphics[width=0.95\linewidth]{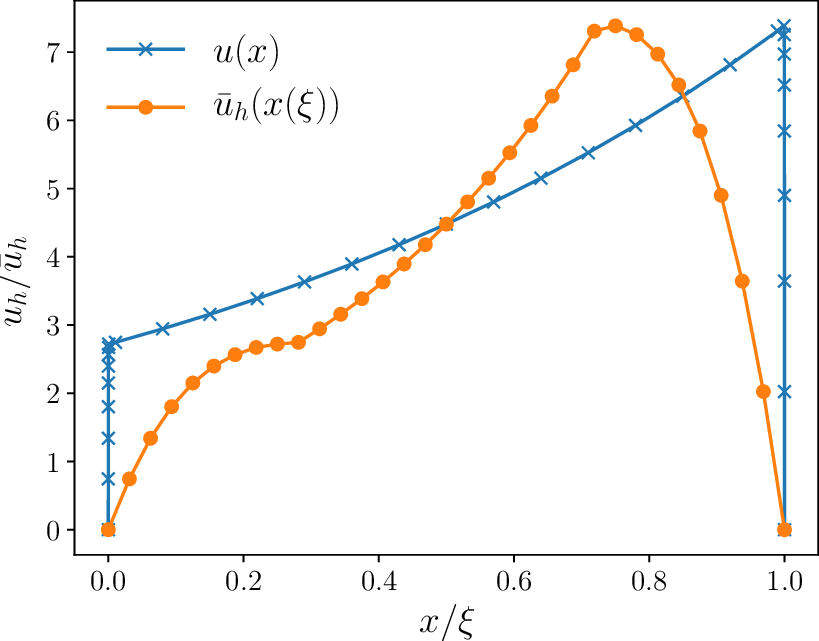}
	\caption{$\varepsilon=10^{-8}$ and $\mu = 10^{-8}$}
	\label{Fig: 1DRCD eps0.0001}
\end{subfigure}
\caption{FEM solutions to~\eqref{eqn: 1DRCD} under each regime with $N=32$ on $\Omega$ and $\Omega^{[c]}$}
\label{Figure: FEM sols RCD}
\end{figure}

We now consider a specific sets of problems of the form
\begin{multline}%
  \label{eqn: 1DRCD}
  -\varepsilon u''(x) + \mu u'(x)+ u(x) = e^{1+x},\\ \text{ for }
  x \in \Omega :=(0,1),\quad\text{with }u(0)=u(1)=0, 
\end{multline}
where $0<\varepsilon\ll 1$ and $0<\mu\leq 1$. This problem  is interesting because 
the  location and width of
layers present in solutions depend on the relative magnitude of $\varepsilon$
and $\mu$, 
leading to three distinct  regimes, as summarised in \autoref{Table:
  regimes} (see~\citep[\S3.2]{Lins2010}). 
Examples of solutions for each regime  on both the physical domain,
$\Omega$, and the computational domain, $\Omega^{[c]}$, are
shown in \autoref{Figure: FEM sols RCD}.
\begin{table}[tbh]
  \centering
  \caption{Location and width of layers, depending on the magnitude of
    $\varepsilon$ and $\mu$.}
  \label{Table: regimes}
  \begin{tabular}{l l c c l}
    \toprule
		Case &Regime&\multicolumn{2}{c}{rate of decay is}&\\
		&& \multicolumn{2}{c}{determined by}& \\
& & $x=0$  & $x=1$ & \\
    \midrule
(a) &    $\varepsilon\ll \mu=1$  & $1$& $1/\varepsilon$ & (convection-diffusion)\\
(b) &  $\varepsilon\ll \mu^2\ll 1$ & $1/\mu$& $\mu/\varepsilon$ &(reaction-convection-diffusion)\\
(c) &  $\mu^2 \ll\varepsilon\ll 1$&  $1/\sqrt{\varepsilon}$ & $1/\sqrt{\varepsilon}$ & (reaction-diffusion)\\
    \bottomrule			
  \end{tabular}
\end{table}

In~\citep{BrZa2016} it is proven that for one-dimensional
reaction-convection-diffusion equations the errors measured in the
energy norm are bounded as 
\begin{equation}%
  \label{eqn: error bound RCD}
  \|u-u_h\|_{\mathrm{E}}^{} \leq
  C\left(\left(\varepsilon^{1/2}+\mu\right)^{1/2}N^{-1} + N^{-2}\right),
\end{equation} 
where $u$ is the exact solution, $u_h$ is the $\mathcal{P}_1$-FEM
solution to the SPDE, computed  on a Bakhvalov-type mesh and
\begin{equation*}
  \|u\|_{\mathrm{E}}^{} := \left(\varepsilon\|u'\|_0^2+\|u\|_0^2\right)^{1/2}.
\end{equation*}
A Bakhvalov mesh for this problem equidistributes \eqref{eq:rho 1DRD}
where $\lambda_0$ and $\lambda_1$ are related to (ordered) solutions
of
\[
  -\varepsilon \lambda(x)^2 + \mu b(x)\lambda(x) + r(x)=0.
\]

To verify that solutions computed on the meshes produced by our
algorithm satisfy the same bounds as in~\eqref{eqn: error bound RCD}, 
we estimate the errors, $\|E_h\|_{\mathrm{E}}^{} :=
\|u_h-u_{2,h}\|_{\mathrm{E}}^{}$,
where   $u_{2,h}$ is  the $\mathcal{P}_2$-FEM
solution to the SPDE, computed on the same mesh as $u_h$.

First, in \autoref{Table: Errors RCD mu=0.001}, we consider the
convergence of the scheme by presenting results over a range of
values of $N$. We have fixed $\mu=10^{-3}$, and taken
various values for $\varepsilon$, so
that each of the three regimes is represented. The errors for~\eqref{eqn: 1DRCD} are
reported, along with the rates of convergence. The number of
iterations of the MPDE 
method performed on the final mesh size is shown in parentheses.
For these calculations, we have taken
$K=0.28$ (resulting in approximately 30\% of mesh points in each layer region) and  $\sigma=2.5$ in~\eqref{eqn: rho dx},
and $\mathrm{TOL}=10^{-3}$ in \Cref{Alg: 1D MPDE dx}.
The results show that the method yields a consistently robust order of convergence and that the number of iterations needed on the final mesh is independent of $\varepsilon$ and $N$.
\begin{table}[htb]
 \centering
  \caption{$\|E_h\|_{\mathbb{E}}^{}$ for \eqref{eqn: 1DRCD} with
    $\mu=10^{-3}$, and rates of convergence (number of iterations)}
  \label{Table: Errors RCD mu=0.001}
  \begin{small}
  \begin{tabular}{l c c c c c c}
    \toprule
    $\varepsilon$& $N=32$ & $N=64$ & $N=128$ & $N=256$ & $N=512$ & $N=1024$\\
    \midrule
1 & 4.08e-02 & 2.04e-02 & 1.02e-02 & 5.11e-03 & 2.55e-03 & 1.28e-03 \\
& (2) & 1.00 (2) & 1.00 (2) & 1.00 (2) & 1.00 (2)& 1.00 (2)\\
$10^{-2}$ & 1.08e-01 & 5.34e-02 & 2.65e-02 & 1.32e-02 & 6.61e-03 & 3.30e-03 \\
& (3) & 1.02 (3) & 1.01 (2) & 1.00 (2)& 1.00 (2)& 1.00 (2)\\
$10^{-4}$ & 5.44e-02 & 2.63e-02 & 1.31e-02 & 6.50e-03 & 3.24e-03 & 1.62e-03 \\
& (6) & 1.05 (3) & 1.01 (3) & 1.01 (2)& 1.00 (2)& 1.00 (2)\\
$10^{-6}$ & 2.27e-02 & 1.07e-02 & 5.30e-03 & 2.66e-03 & 1.33e-03 & 6.65e-04 \\
& (8) & 1.08 (3)& 1.02 (3)& 1.00 (2)& 1.00 (2)& 1.00 (2)\\
$10^{-8}$ & 1.82e-02 & 8.35e-03 & 4.10e-03 & 2.04e-03 & 1.02e-03 & 5.12e-04 \\
& (9)  & 1.13 (4) & 1.03 (3)& 1.00 (2)& 1.00 (2) & 1.00 (2)\\
$10^{-10}$ & 1.80e-02 & 8.33e-03 & 4.07e-03 & 2.04e-03 & 1.02e-03 & 5.10e-04 \\
& (9) & 1.12 (4)& 1.03 (3)& 1.00 (3)& 1.00 (3)& 1.00 (2)\\
$10^{-12}$ & 1.80e-02 & 8.36e-03 & 4.08e-03 & 2.04e-03 & 1.02e-03 & 5.10e-04 \\
&  (9)& 1.10 (5)& 1.03 (3)& 1.00 (3)& 1.00 (3)& 1.00 (3)\\

    \bottomrule
	\end{tabular}
      \end{small}
    \end{table}

Since the results in a given column of \autoref{Table: Errors RCD
  mu=0.001} straddle multiple regimes, the robustness of the error
estimates, with respect to the perturbation parameters,
might not be so clear. Therefore,  in \autoref{Table: vary
  mu}, we present results for a single value of $N$. Moreover,  each
row is restricted to a single case (as per \autoref{Table: regimes})
by  taking a fixed value of $\mu$,
and only values of $\varepsilon$ that correspond to the associated regime.
For Cases (a) and (b), we see that the error is clearly robust with
respect to $\varepsilon$, (suggesting the $\mu$-weighted term in~\eqref{eqn: error bound RCD}
  dominates). For Case (c), the error
scales as $\varepsilon^{1/4}$, as expected. We can conclude that the error
bound~\eqref{eqn: error bound RCD} is satisfied.
\begin{table}[htb]
  \centering
  \caption{$\|E_h\|_{\mathrm{E}}^{}$ for $N=1024$ for each case in \autoref{Table: regimes}}
  \label{Table: vary mu}
    \begin{footnotesize}
  \begin{tabular}{lccccccc}
    \toprule
    Case~(a)& $\varepsilon=10^{-4}$& $\varepsilon=10^{-5}$ & $\varepsilon=10^{-6}$    & $\varepsilon=10^{-7}$ & $\varepsilon=10^{-8}$&$\varepsilon=10^{-9}$&$\varepsilon=10^{-10}$\\
      $\mu=1$	&9.78e-03 & 9.73e-03 & 9.77e-03 & 9.78e-03 & 9.68e-03 & 9.68e-03 & 9.64e-03 \\
    \midrule 
        Case~(b) & ---&---& $\varepsilon=10^{-8}$ & $\varepsilon=10^{-9}$& $\varepsilon=10^{-10}$& $\varepsilon=10^{-11}$&$\varepsilon=10^{-12}$\\
     $\mu=10^{-3}$& && 5.12e-04 & 5.10e-04 & 5.10e-04 & 5.10e-04 & 5.10e-04 \\     
     \midrule
    Case~(c) &---&---&---& $\varepsilon=10^{-6}$ & $\varepsilon=10^{-8}$&$\varepsilon=10^{-10}$&$\varepsilon=10^{-12}$\\ 
     $\mu=10^{-8}$ &&&& 5.44e-04 & 1.71e-04 & 5.40e-05 & 1.70e-05 \\
     
		\bottomrule
  \end{tabular}
\end{footnotesize}
\end{table}

\section{Two-dimensional problems}%
\label{sec: 2D problem}
We now extend the approach presented in \S\ref{sec: 1D problems} to
problems in two-dimensions. Again, at the  core of the idea is the
technique for generating layer-adapted meshes, of the Bakhvalov type,
using MPDEs in~\citep{HiMa2021}, but now using iteratively computed
derivative estimates to automatically detect the location and width of layers.

\subsection{Two-dimensional SPDEs}\label{sec: 2D SPDE}
We consider the family two-dimensional reaction-convection-diffusion
equations of the form
\begin{multline}
  \label{eqn: 2DRCD}
  -\varepsilon\Delta u(x,y) + \mu\mathbf{b}\cdot\nabla u(x,y) + u(x,y) = f(x,y),\\\text{ for } (x,y)\in\Omega:=(0,1)^2 \text{ and } u\rvert_{\partial \Omega}=0.
\end{multline}
This motivates the two-dimensional MPDE formulation discussed in
\S\ref{sec: 2D MPDE dx}. In particular, the location and width of the layers in the solutions to~\eqref{eqn:
  2DRCD} is dependent on the relative magnitude of $\varepsilon$ and
$\mu$, as outlined for the one-dimensional
problem in \autoref{Table: regimes}.
However, there are two further complications: layers may occur at any
of the four boundaries, and, depending on the direction of
$\mathbf{b}$, boundary layers maybe exponential or parabolic in
nature.

\begin{figure}[tbh]
	\centering
	\begin{subfigure}{0.5\linewidth}
	\centering 	
	\includegraphics[width = 0.95\linewidth]{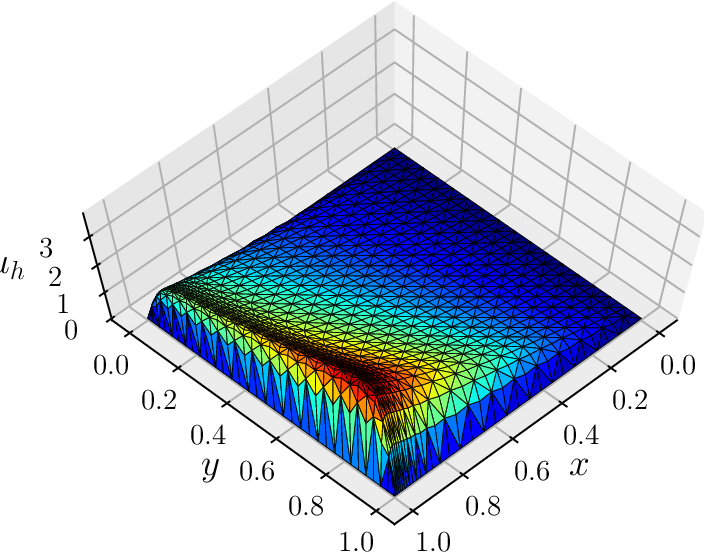}
	\caption{$\varepsilon =10^{-2}$ and $\mu=1$}
	\label{Fig: eps 1e-2, mu 1}
\end{subfigure}  
	\begin{subfigure}{0.48\linewidth}
	\centering 	 	
	\includegraphics[width = 0.95\linewidth]{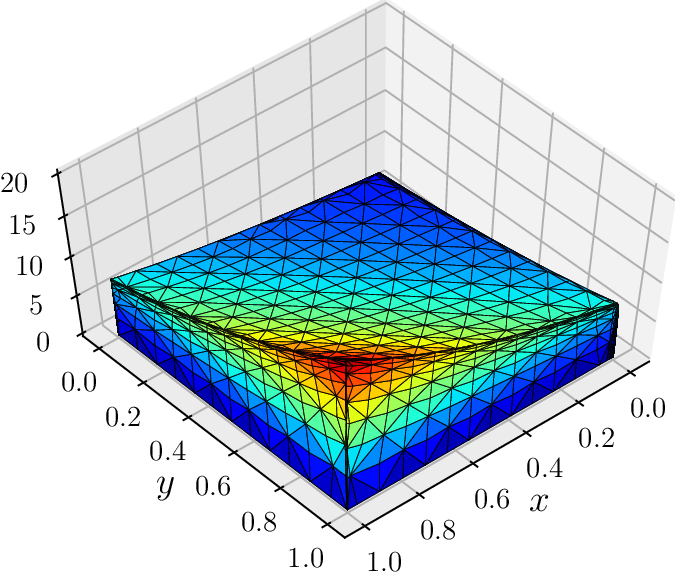}
	\caption{$\varepsilon =10^{-8}$ and $\mu=10^{-3}$}
	\label{Fig: eps 1e-8, mu 1e-3}
\end{subfigure}  
	\begin{subfigure}{0.48\linewidth}
		\centering
		\includegraphics[width = 0.95\linewidth]{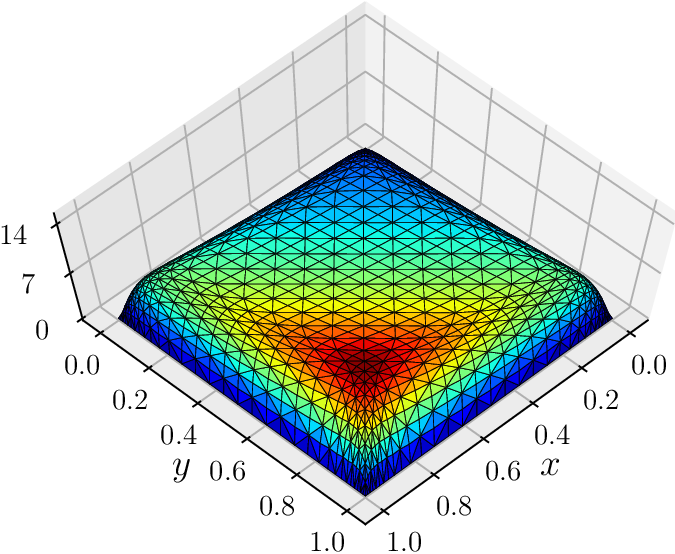}
		\caption{$\varepsilon =10^{-3}$ and $\mu=10^{-3}$}
		\label{Fig: eps 1e-3, mu 3}
	\end{subfigure}
	\caption{Examples of the FEM solution for~\eqref{eqn: 2DRCD} with $N=32$}
	\label{Fig: Solution 2DRCD}
\end{figure}

For a specific example, we consider the  problem
\begin{multline}
  \label{eqn: 2DRCD example}
  -\varepsilon\Delta u(x,y) + \mu(3-x)u_x + u(x,y) = e^{1+x+y},\\\text{ for } (x,y)\in\Omega:=(0,1)^2 \text{ and } u\rvert_{\partial \Omega}=0.
\end{multline}
The possible variation is the solutions is demonstrated in
\autoref{Fig: Solution 2DRCD}  which presents 
solutions to \eqref{eqn: 2DRCD} for a selection of values for
$\varepsilon$ and $\mu$. Note that the relative magnitude of $\varepsilon$
and $\mu$ impacts  the number of layers that are visible, their location,
and their width.

\subsection{Two-dimensional MPDEs}
\label{sec: 2D MPDE dx}
For~\eqref{eqn: 2DRCD}, a mesh that may be
layer-adapted at any of the four boundaries is needed and so the
MPDE is formulated as the two-dimension vector-valued Poisson equation, for $\boldsymbol{x}(\xi_1, \xi_2) = (x,y)^T$,
\begin{subequations}
  \label{eqn: 2D MPDE generic dx}
\begin{equation}%
  -\nabla^{[c]}\cdot(M(\boldsymbol{x}(\xi_1,\xi_2), u_h)\nabla^{[c]}
  \boldsymbol{x}(\xi_1,\xi_2)) = (0,0)^T,\text{ for } (\xi_1,
  \xi_2) \in \Omega^{[c]}, 
\end{equation}
subject to the  boundary conditions
\begin{equation}%
	\label{eqn: 2D MPDE bcs dx}
	\begin{split}
		x(0,\xi_2) = 0, \quad 
		x(1,\xi_2) = 1, \quad 
		\frac{\partial x}{\partial \mathbf{n}}(\xi_1,0) = 0, \quad 
		\frac{\partial x}{\partial \mathbf{n}}(\xi_1,1) = 0, \\
		y(\xi_1,0) = 0, \quad 
		y(\xi_1,1) = 1, \quad 
		\frac{\partial y}{\partial  \mathbf{n}}(0,\xi_2) = 0, \quad 
		\frac{\partial y}{\partial  \mathbf{n}}(1,\xi_2) = 0. \\
	\end{split}
\end{equation}  
\end{subequations}
As before, the parameters  $K_1$ and $K_2$ determine the proportion
of mesh points in the layer regions, and $\sigma$ is based on the
formal order of them scheme.

\begin{figure}[tbh]
  \centering
  \begin{subfigure}{0.45\linewidth}	
    \includegraphics[width=0.89\linewidth]{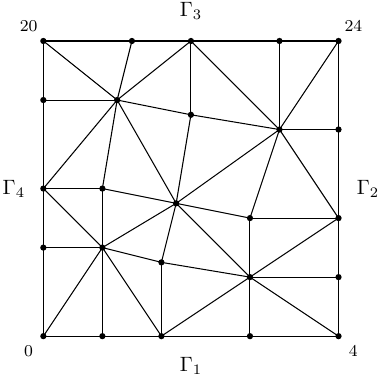}
    \caption{Boundaries}
    \label{Fig: mesh boundaries}
  \end{subfigure}
  \begin{subfigure}{0.48\linewidth}	
    \includegraphics[width=0.9\linewidth]{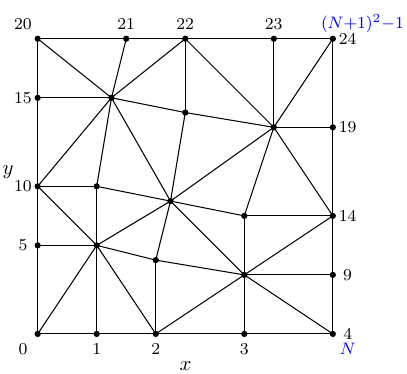}
    \caption{Indices of the mesh points}
    \label{Fig: mesh indexing}
  \end{subfigure}
  \caption{Example of boundaries and indices of mesh points on a non-tensor product grid with $N=4$}
  \label{Fig: 5*5 mesh}
\end{figure}

The computational meshes that are to be computed are not necessarily
tensor-product, however they are $M$-uniform. Therefore, we can denote
the mesh points as  
$(x, y)_k$, for $k = 0,1,\dots, (N+1)^2-1$, indexed  using standard
lexicographic ordering. For example, in \autoref{Fig: mesh indexing} we show the indices of the mesh points when $N=4$.

In \eqref{eqn: 2D MPDE generic dx}, we take
\begin{subequations}\label{eq:M etc}
\begin{equation}%
  \label{eqn: M a posteriori}
  M(\boldsymbol{x}, u_h) =
  \begin{pmatrix}\displaystyle
    m_{1,1} & 0\\
    0 &m_{2,2}
  \end{pmatrix},
\end{equation}
where
\begin{equation}\label{eq:mij}
  \begin{split}
    m_{1,1}(\boldsymbol{x}, u_h) &= \max\Bigg{\{}1,K_1\left(\upsilon_4\exp\left(-\upsilon_4 x\sigma^{-1}\right)+\upsilon_2 \exp\left(-\upsilon_2(1-x)\sigma^{-1}\right)\right)\Bigg{\}}\\
    m_{2,2}(\boldsymbol{x}, u_h) &=
      \max\Bigg{\{}1,K_2\left(\upsilon_1\exp\left(-\upsilon_1 y\sigma^{-1}\right)+\upsilon_3 \exp\left(-\upsilon_3(1-y)\sigma^{-1}\right)\right)\Bigg{\}}.
	\end{split}
      \end{equation}
\end{subequations}
Here $m_{1,1}$ and $m_{2,2}$ can be thought of as extensions to
the terms in \eqref{eqn: rho dx} to two dimension: they encode that
there may be layers adjacent the each of four boundaries. The terms
$v_1$, $v_2$, $v_3$ and $v_4$ determine the magnitude and decay rate
of the associated layer terms. Taken together, they relate to
pointwise bounds on derivatives of $u(x,y)$~\citep{HiMa2021}. Since
each layer term decays rapidly away from its associated boundary, and
takes its maximum at that boundary, the $v_i$ terms are determined by
the appropriate derivatives of the numerical solution adjusted for the
effects of the right-hand side of the SPDE, evaluated at the
boundary. For efficiency, these are computed as  
\begin{subequations}%
  \label{eqn: v_xy} 
  \begin{align}
    \upsilon_1(\boldsymbol{x}, u_h)_k &=\frac{|(u_h)_{y^+}(x, y)_{l}|}{\max\{1,|f(x,y)_l|\}}, \text{ where } l = k  \text{ mod } N{+}1,\\
    \upsilon_2(\boldsymbol{x}, u_h)_k &=\frac{|(u_h)_{x^-}(x, y)_{l}|}{\max\{1,|f(x,y)_l|\}}, \text{ where } l= N + \bigg{\lfloor}\frac{k}{N{+}1}\bigg{\rfloor}(N{+}1),\\
    \upsilon_3(\boldsymbol{x}, u_h)_k &=\frac{|(u_h)_{y^-}(x, y)_{l}|}{\max\{1,|f(x,y)_l|\}}, \text{ where } l = N(N{+}1)+k  \text{ mod } N{+}1,\\
    \upsilon_4(\boldsymbol{x}, u_h)_k &=\frac{|(u_h)_{x^+}(x, y)_{l}|}{\max\{1,|f(x,y)_l|\}}, \text{ where } l = \bigg{\lfloor} \frac{k}{N+1}\bigg{\rfloor}(N+1),	
  \end{align}
\end{subequations}
where $k = 0,1,\dots, (N+1)^2-1$, are the (lexicographic) indices of
the mesh points,
and $f(x,y)$ is the right-hand side of the SPDE.
That is, these functions propagate the derivatives of $u_h$, adjusted by
the value of $f(x,y)$ at their respective boundaries across the
domain. 
For example,  in the case where $N=4$, we take
$\upsilon_4(\boldsymbol{x},u_h)$ to have the same value at each of the mesh
points $k=5$, $k=6$, \dots, $k=9$ (see \autoref{Fig: mesh indexing}),
which is $|(u_h)_{x^+}(x, y)_5|/\max\{1,|f(x,y)_5|\}$.
If one preferred, we could set, for example,
$\upsilon_4(\boldsymbol{x},u_h) = |(u_h)_{x^+}(0,y)|/\max\{1,|f(0,y)|\}$. However, we have verified that there is no noticeable
advantage for the extra computational expense one would encounter on 
a non-tensor product grid.

Finally, we note that if a particular edge does not feature a layer
the  associate $v_i$ term will be $\mathcal{O}(1)$, and so does not induce any refinement near that boundary.

\subsection{Algorithm and implementation}%
\label{sec: 2D algorithm dx V2}
Although the approach to numerically solving~\eqref{eqn: 2D MPDE
  generic dx} efficiently is an extension of that used for the
one-dimensional problem in \S\ref{sec: Algorithm for MMPDE dx}, it
has minor variations that  merit further discussion. (For completeness
the full algorithm is given in \autoref{app: 2D algorithm}).
As in~\Cref{Alg: 1D MPDE dx}, we take an initial mesh with $N=16$.
We then preform  $\big{\lceil}
6\log_{10}\left(\varepsilon^{-1}\right)\big{\rceil}$ iterations of
the fixed-point method to resolve the MPDE sufficiently. 
Since only one mesh point is added to a layer region at each
iteration, and we require the minimum mesh width to be
$\mathcal{O}(\varepsilon/N )$, one can deduce that $\mathcal{O}(\log
\varepsilon^{-1})$ iterations are required.
Experiments indicate that $\lceil 6 \log 10 \varepsilon^{-1}\rceil $
iterations suffice.

We then alternate between a sequence of uniform $h$-refinements of the
computational mesh, and 5 iterations of the fixed point method for the
MPDE on each of these, until the mesh has the required number of mesh
points.

Computational experimentation indicated that this approach is
sufficient, and more efficient than
iterating until a specified tolerance is achieved.

\subsection{Numerical results}%
\label{sec: 2D numerical results}
The mesh for~\eqref{eqn: 2DRCD example} is generated using~\Cref{Alg:
  2D MPDE a posteriori} with MPDE~\eqref{eqn: 2D MPDE generic dx}. We
set $K_1= K_2 = 0.28$ which results in
approximately $30\%$ of the mesh points being located in each layer
region.
One takes $\sigma>p+1$, where $p$ is the order of the FEM, so we have
set $\sigma=2.5$ for  this $\mathcal{P}_1$-FEM.

\begin{figure}[tbh]
	\centering
	\begin{subfigure}{0.42\linewidth}
		\includegraphics[height=0.95\linewidth]{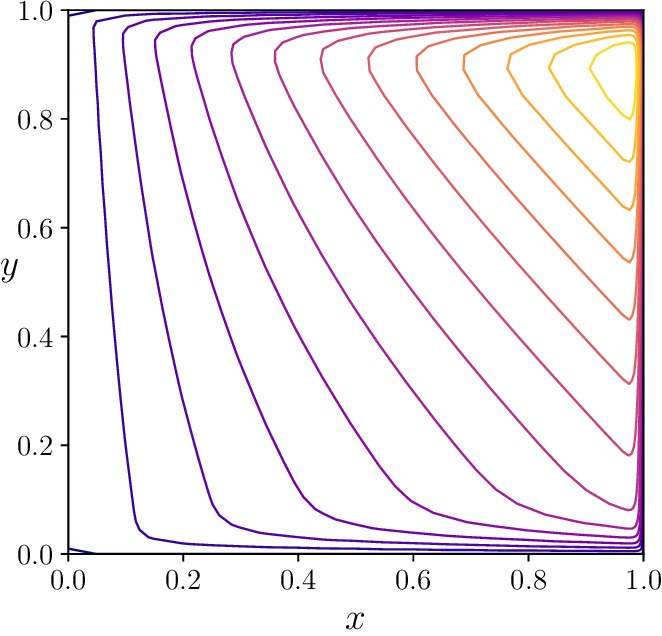}
		\caption{Contour plot on $\Omega$}
		\label{Fig: contour plot omega}
	\end{subfigure}
	\begin{subfigure}{0.42\linewidth}
		\includegraphics[height=0.95\linewidth]{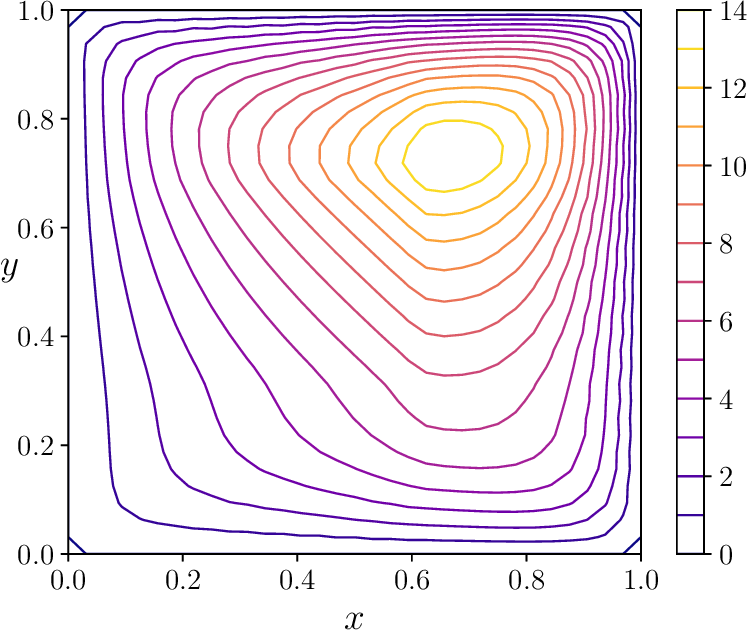}
		\caption{Contour plot on $\Omega^{[c]}$}
		\label{Fig: contour plot omega^c}
	\end{subfigure}
	\caption{Contour plots with $N=32$, $\varepsilon = 10^{-3}$, and $\mu = 10^{-1}$}
	\label{Fig: contour plots}
\end{figure}

An example of the mesh for~\eqref{eqn: 2DRCD example} generated using ~\eqref{eqn: 2D MPDE generic dx}  is shown in \autoref{Fig:
  mesh2D} and one observes that this is not a tensor-product grid. In
\autoref{Fig: contour plot omega^c}, one sees that the layer regions
in the related SPDE are resolved when transformed onto the
computational domain, $\Omega^{[c]}$. 
\begin{figure}[tbh]
	\centering
	\includegraphics[width=0.5\linewidth]{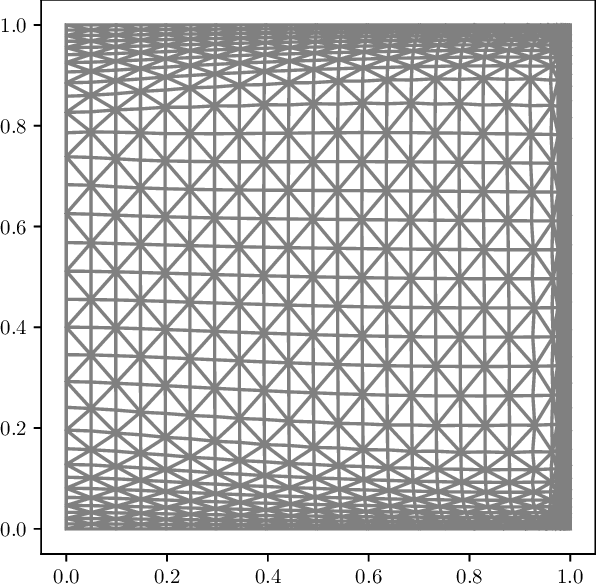}
	\caption{Mesh for~\eqref{eqn: 2DRCD example} with $N=32$, $\varepsilon=10^{-3}$ and $\mu=10^{-1}$}
	\label{Fig: mesh2D}	
\end{figure}

As in \S\ref{sec:1D numerics}, we first  fix $\mu=10^{-3}$, and
take various values for $\varepsilon$, so that each of the three regimes is
represented. The errors measured in the energy norm for~\eqref{eqn:
  2DRCD example} when solved on a mesh generated using~\Cref{Alg: 2D
  MPDE a posteriori} are shown in \autoref{Table: Errors 2DRCD}.
These clearly show that results are robust with respect to $\varepsilon$,
and first-order  convergent  with respect to $N$
\begin{table}[tbh]
	\centering
	\caption{$\|E_h\|_{\mathrm{E}}^{}$ for~\eqref{eqn: 2DRCD example}
		solved on meshes generated using~\Cref{Alg: 2D MPDE a posteriori},
		with $\mu=10^{-3}$, and various values of $N$ and $\varepsilon$}
	\label{Table: Errors 2DRCD}
	\begin{small}
		\begin{tabular}{l c c c c c }
			\toprule
			$\varepsilon /N$ &32&64&128&256&512\\
			\midrule		
1 & 8.19e-02 & 4.10e-02 & 2.05e-02 & 1.03e-02 & 5.13e-03 \\
&   & 1.00 & 1.00 & 1.00 & 1.00 \\
$10^{-2}$ & 2.79e-01 & 1.38e-01 & 6.83e-02 & 3.40e-02 & 1.69e-02 \\
&    & 1.02 & 1.01 & 1.01 & 1.00 \\
$10^{-4}$ & 1.22e-01 & 5.97e-02 & 2.96e-02 & 1.47e-02 & 7.34e-03 \\
&    & 1.03 & 1.01 & 1.01 & 1.00 \\
$10^{-6}$ & 5.14e-02 & 2.49e-02 & 1.23e-02 & 6.13e-03 & 3.06e-03 \\
&    & 1.05 & 1.01 & 1.01 & 1.00 \\
$10^{-8}$ & 4.03e-02 & 1.94e-02 & 9.66e-03 & 4.80e-03 & 2.39e-03 \\
&    & 1.05 & 1.01 & 1.01 & 1.01 \\
$10^{-10}$ & 3.92e-02 & 1.90e-02 & 9.45e-03 & 4.70e-03 & 2.34e-03 \\
&    & 1.04 & 1.01 & 1.01 & 1.00 \\
$10^{-12}$ & 3.91e-02 & 1.90e-02 & 9.44e-03 & 4.69e-03 & 2.34e-03 \\
&    & 1.04 & 1.01 & 1.01 & 1.01 \\
			\bottomrule
		\end{tabular}
	\end{small}
\end{table}

In \autoref{Table: vary mu 2D}, we fix $N=512$, and examine the
results for various values of $\varepsilon$ and $\mu$. Again, we see robust
convergence  with $\varepsilon$ in each different regime. It should be noted
that, for $\mu=1$ and $\mu=10^{-3}$, the problem is essentially
convection-dominated, and so there is no $\varepsilon$ dependency in the
computed errors. For $\mu=10^{-8}$, the problem is dominated by the
reaction term, and so the $\varepsilon$-dependency in the error is consistent
with \eqref{eqn: error bound RCD}, and with those
shown in \autoref{Table: vary mu} for the one-dimensional problem. /

\begin{table}[H]
	\centering
	\caption{$\|E_h\|_{\mathrm{E}}^{}$ for~\eqref{eqn: 2DRCD example} with $N=512$ for each case in \autoref{Table: regimes}}
	\label{Table: vary mu 2D}
	\begin{footnotesize}
		\begin{tabular}{lccccccc}
			\toprule
Case (a) & $\varepsilon = 10^{-4}$ & $\varepsilon = 10^{-5}$ & $\varepsilon = 10^{-6}$ & $\varepsilon = 10^{-7}$ & $\varepsilon = 10^{-8}$ & $\varepsilon = 10^{-9}$ & $\varepsilon = 10^{-10}$ \\
$\mu = 1$ & 3.90e-02 & 3.49e-02 & 3.40e-02 & 3.39e-02 & 3.28e-02 & 3.46e-02 & 3.29e-02 \\
\midrule
Case (b) &  &  & $\varepsilon = 10^{-8}$ & $\varepsilon = 10^{-9}$ & $\varepsilon = 10^{-10}$ & $\varepsilon = 10^{-11}$ & $\varepsilon = 10^{-12}$ \\
$\mu=10^{-3}$ &  &  & 2.39e-03 & 2.35e-03 & 2.34e-03 & 2.34e-03 & 2.34e-03 \\
\midrule
Case (c) &  &  &  & $\varepsilon = 10^{-6}$ & $\varepsilon = 10^{-8}$ & $\varepsilon = 10^{-10}$ & $\varepsilon = 10^{-12}$ \\
$\mu = 10^{-8}$ &  &  &  & 2.44e-03 & 7.84e-04 & 2.49e-04 & 7.91e-05 \\
			\bottomrule
		\end{tabular}
	\end{footnotesize}
\end{table}

\section{Conclusions and future work}%
\label{sec: conclusion MPDE dx}	
The MPDEs presented in
\eqref{eqn: 1DMPDE dx} and \eqref{eqn: 2D MPDE generic dx} use \emph{a
  posteriori} information about the related SPDE. More precisely,
using only the knowledge that the solution possesses boundary layers
(but not their location or width) layer-adapted meshes are
generated. The solutions are robust and the errors converge as
expected. The magnitude and rates of convergence of the errors are
similar to when the solutions are generated using the MPDE method
based on \emph{a priori} information~\citep{HiMa2021}. Our
investigations included ensuring that the local mesh width is
appropriate for the relevant layer width, though, for brevity, we have
not included the detail.

There are numerous possibilities for extending this work. Perhaps the
most obvious, and challenging, is to generalise the approach to
produce layer-adapted meshes for interior layer problems. Work in this
direction is currently in its infancy.

\medskip

\noindent \textbf{Acknowledgement:} The work of RH is supported by the Irish
Research Council, GOIPG/2017/463 \&  GOIPD/2022/284.

\clearpage
\begin{appendices}
  \section{FEniCS code to compute a 1D layer-adapted
    mesh}\label{app: code 1D}
  Here we present Python code an implementation of \autoref{Alg: 1D MPDE dx}
  for solving 
\begin{equation*}%
  -\varepsilon u''(x) + \mu b(x) u'(x) + r(x) u(x) =  f(x)\quad \text{for }
  x \in (0,1), \quad\text{ and } u(0)=u(1)=0,
\end{equation*}
  for the specific example in \eqref{eqn: 1DRCD} with $\varepsilon=10^{-3}$, $\mu=10^{-3}$ and $N=32$.
We use the Gauss Lobatto quadrature rule to solve both the MPDE and
SPDE. (See the note in ~\citep[App.~B]{HiMa2021} for why this is necessary).
\begin{lstlisting}
## 1DMPDE_dx_RCD_code.py
# FEniCS version: 2019.2.0.dev0
# Generate a layer-adapted mesh for the one-dimensional
# reaction-convection-diffusion equation
#    -eps u'' + mu u' + u = exp(1+x) for x in (0,1) with u(0)=u(1)=0,
# via mesh partial differential equations and a posteriori information using
#   -(rho(x) x'(xi))' = 0 for xi in (0,1) with x(0) = x(1) = 1.
from fenics import *
import math; import numpy as np
# Change quadrature rule to Gauss Lobatto
from FIAT.reference_element import *
from FIAT.quadrature import *
from FIAT.quadrature_schemes import create_quadrature
def create_GaussLobatto_quadrature(ref_el, degree, scheme="default"):
    if degree < 3: degree = 3
    return GaussLobattoLegendreQuadratureLineRule(ref_el, degree)
import FIAT
FIAT.create_quadrature = create_GaussLobatto_quadrature

# Parameters for physical problem and mesh size
epsilon = 1.0e-3                          # coefficient of diffusion term
mu = 1.0e-3                               # coefficient of convection term
f = Expression('exp(1.0+x[0])', degree=2) # RHS of SPDE
N = 2**5                                  # Max mesh intervals; at least 2**4

# Parameters for MPDE
K = 0.28                            # proportion of mesh points in layers
sigma = 2.5                               # related to degree of elements
x_init = Expression('x[0]', degree=2)     # Initial solution to the MPDE

# Parameters for stepping through N
N_start = 2**4                            # initial mesh size
TOL = 1.0e-3                              # Stopping criterion tolerance

# Lists to store meshes and function spaces (as we step through mesh sizes)
meshc_list, mesh_list = [None]*int(math.log(N,2)-3),[]
Vc_list, V_list = [None]*int(math.log(N,2)-3),[]

def comp_space(N): # Computational function space parameters
    meshc = UnitIntervalMesh(N)
    Vc = FunctionSpace(meshc, 'P', 1)     # Function space with P1-elements
    vc = TestFunction(Vc)                 # Test function on Vc
    xc = TrialFunction(Vc)                # Trial function on Vc
    return meshc, Vc, vc, xc

def MPDE_FEM(rho, x_old, u_h_x0, u_h_xN, xc, vc, meshc, Vc): # Solve MPDE
    bc = DirichletBC(Vc, x_init, 'on_boundary')   
    a = rho(x_old, u_h_x0, u_h_xN)*inner(grad(xc),grad(vc))*dx(meshc) # LHS of MPDE
    L = Constant("0.0")*vc*dx             # MPDE is homogeneous
    x_new = Function(Vc)                  # New mesh function
    solve(a==L, x_new, bc)
    return x_new

def rho(x, u_h_x0, u_h_xN): # Define rho(x) for MPDE
    upsilon_0 = u_h_x0/max(1.0,f(0.0))
    upsilon_1 = u_h_xN/max(1.0,f(1.0))
    return conditional((K*(upsilon_0*exp(-upsilon_0*x/sigma)+upsilon_1*exp(-upsilon_1*(1-x)/sigma)))>1.0,\
        (K*(upsilon_0*exp(-upsilon_0*x/sigma)+upsilon_1*exp(-upsilon_1*(1-x)/sigma))), 1.0)

def phys_space(N, x): # Physical function space parameters
    mesh = UnitIntervalMesh(N)
    mesh.coordinates()[:,0] = x.compute_vertex_values()[:]
    V = FunctionSpace(mesh, 'P', 1)        # Function space with P1-elements
    v = TestFunction(V)                    # Test function on V
    u = TrialFunction(V)                   # Trial function on V
    return mesh, V, v, u

def SPDE_FEM(u, v,  mesh, V): # Solve reaction-convection-diffusion problem
    bcp = DirichletBC(V, "0.0", 'on_boundary') 
    a = epsilon*dot(grad(u), grad(v))*dx(mesh) \
      + mu*u.dx(0)*v*dx + u*v*dx
    L = f*v*dx
    u_h = Function(V)
    solve(a==L, u_h, bcp)
    # calculate the derivatives of u_h on the first and last intervals
    u_h_x = (project(u_h.dx(0), V)).compute_vertex_values()
    u_h_x0, u_h_xN = abs(u_h_x[0]), abs(u_h_x[-1])
    return u_h, u_h_x0, u_h_xN

# Initial values
u_h_x0, u_h_xN = 1.0, 1.0
meshc_list[0], Vc_list[0], vc, xc = comp_space(N_start)
x = interpolate(x_init, Vc_list[0])        # Initial value for x

k = 0
for i in range(0,int(math.log(N,2)-3)): # Iterate through h-refinements
    N_step = 2**(i+4)                      # mesh size
    Delta = TOL+1
    while Delta > TOL:
        mesh_list.append(None), V_list.append(None)
        mesh_list[k], V_list[k], v, u = phys_space(N_step, x)
        u_h_x0_old, u_h_xN_old = u_h_x0, u_h_xN
        u_h, u_h_x0, u_h_xN = SPDE_FEM(u, v,  mesh_list[k], V_list[k])
        Delta = abs((u_h_x0-u_h_x0_old)/u_h_x0)+abs((u_h_xN-u_h_xN_old)/u_h_xN)
        x = MPDE_FEM(rho, x, u_h_x0, u_h_xN, xc, vc,meshc_list[i],Vc_list[i])
        k = k+1
    if N_step < N:
        # Generate computational function space on the finer mesh
        meshc_list[i+1], Vc_list[i+1], vc, xc = comp_space(N_step*2)
        # Interpolate the solution onto the new computational function space
        x = interpolate(x, Vc_list[i+1])
        # Generate physical function space on the finer mesh
        mesh_list[i+1], V_list[i+1], v, u = phys_space(N_step*2, x)

# final (adapted) mesh and solution to SPDE
mesh_list.append(None), V_list.append(None)
mesh_list[k], V_list[k], v, u = phys_space(N, x)
u_h, u_h_x0, u_h_xN = SPDE_FEM(u, v,  mesh_list[k], V_list[k])
\end{lstlisting}
\clearpage

  \section{2D Algorithm}%
  \label{app: 2D algorithm}
\begin{algorithm}[H]%
  \SetKwInOut{Inputs}{Inputs}
  \caption{Generate a two-dimensional layer-adapted mesh using an MPDE using \emph{a posteriori} information and $h$-refinement.}
  \label{Alg: 2D MPDE a posteriori}
  \Inputs{$N$, the number of intervals in both directions in the final mesh on $\Omega^{[c]}$; \textit{M}, a mesh density function.} 
  Set $\omega_h^{[c;0]} := \{(\xi_1, \xi_2)_i\}_{i=0}^{16}$ to be a uniform tensor-product mesh on $\overline{\Omega}^{[c]}$ with $16$ mesh intervals in each coordinate direction\;
  Set $\textbf{\textit{x}}(\xi_1,\xi_2) = (\xi_1,\xi_2)$ for $(\xi_1,\xi_2) \in \overline{\Omega}^{[c]}$\;
  $\textbf{\textit{r}} \leftarrow \boldsymbol{x}$;  $k\leftarrow$ 0\;
  \For{$j$ in $0:\big{\lceil} 6\log_{10}\left(\varepsilon^{-1}\right)\big{\rceil}$}{	
    Set $\omega_h^{[k]} := \boldsymbol{x}(\omega_h^{[c;0]})$ to be the adapted mesh on $\overline{\Omega}$ at step $k$\;
    $u_h^{[k]} \leftarrow \mathcal{P}_1$-FEM solution to the SPDE on $\omega_h^{[k]}$\;
    Calculate $\upsilon_l^{[k]}(\boldsymbol{r}, u_h^{[k]})$ for $l =
    1,2,3,4$, and then $M(\textbf{\textit{r}},u_h^{[k]})$, as in \eqref{eq:M etc}\;		
    Set $\boldsymbol{x}$ to be the $\mathcal{P}_1$-FEM
    solution, on $\omega_h^{[c;0]}$, to
    \begin{equation}\label{eqn: 2DMPDE_step1}
      \nabla\cdot(M(\textbf{\textit{r}},u_h^{[k]})\nabla \textbf{\textit{x}}(\xi_1,\xi_2)) = (0,0)^T,\text{ for } (\xi_1,\xi_2) \in \Omega^{[c]},
    \end{equation}
    with boundary conditions as defined in~\eqref{eqn: 2D MPDE bcs dx}\;
    $\boldsymbol{r} \leftarrow \boldsymbol{x}$;
    $k \leftarrow k{+}1$\;
  }
  \For{$i$ in $1:(\log_2(N)-4)$}{
    $w^{[c;i]}\leftarrow $ uniform \emph{h}-refinement of $w^{[c;i{-}1]}$\;
    $\boldsymbol{r} \leftarrow \boldsymbol{r}$ interpolated onto $w^{[c;i]}$\;		
    \For{$j$ in $1:4$}{
      Set $\omega_h^{[k]} := \boldsymbol{x}(\omega_h^{[c;i]})$ to be the adapted mesh on $\overline{\Omega}$ at step $k$\;
      $u_h^{[k]} \leftarrow \mathcal{P}_1$-FEM solution to the SPDE on $\omega_h^{[k]}$\;

      Calculate $\upsilon_l^{[k]}(\boldsymbol{r}, u_h^{[k]})$ for $l
      = 1,2,3,4$ and then $M(\textbf{\textit{r}},u_h^{[k]})$, as in \eqref{eq:M etc}\;
      Set \textbf{\textit{x}} to be the $\mathcal{P}_1$-FEM
      solution, on $\omega_h^{[c;i]}$, to
      \begin{equation}\label{eqn: 2DMPDE_step2}
        \nabla\cdot(M(\textbf{\textit{r}},u_h^{[k]})\nabla \textbf{\textit{x}}(\xi_1,\xi_2)) = (0,0)^T,\text{ for } (\xi_1,\xi_2) \in \Omega^{[c]},
      \end{equation}
      with boundary conditions as defined in~\eqref{eqn: 2D MPDE bcs dx}\;
      $\boldsymbol{r} \leftarrow \boldsymbol{x}$; $k \leftarrow k{+}1$\;			
    }}
  Set $\omega_h := \boldsymbol{x}(\omega_h^{[c;i]})$ to be the adapted mesh on $\overline{\Omega}$\;
\end{algorithm}

\end{appendices}
	

\end{document}